\documentclass[11pt]{article}
\usepackage{amssymb,amsfonts,amsmath,latexsym,epsf,tikz,url}

\newtheorem{theorem}{Theorem}[section]

\newtheorem{corollary}[theorem]{Corollary}
\newtheorem{lemma}[theorem]{Lemma}

\newtheorem{definition}[theorem]{Definition}
\newtheorem{example}[theorem]{Example}
\newcommand{\proof}{\noindent{\bf Proof. }}
\newcommand{\qed}{\hfill $\square$\medskip}

\begin{document}

\title{  On the number of fair dominating sets of graphs }

\author{
Saeid Alikhani$^{}$\footnote{Corresponding author}
\and
Maryam Safazadeh
}

\date{\today}

\maketitle

\begin{center}
Department of Mathematics, Yazd University, 89195-741, Yazd, Iran\\
{\tt alikhani@yazd.ac.ir, msafazadeh92@gmail.com}
\end{center}


\begin{abstract}

Let $G=(V,E)$ be a simple graph. A dominating set of $G$ is a subset $D\subseteq V$ such that every vertex not in $D$ is adjacent to at least one vertex in $D$.
 The cardinality of a smallest dominating set of $G$, denoted by $\gamma(G)$, is the domination number of $G$.
  For $k \geq 1$, a $k$-fair dominating set ($kFD$-set) in $G$, is a dominating set $S$ such that $|N(v) \cap D|=k$ for every vertex $ v \in V\setminus D$.
 A fair dominating set, in $G$ is a $kFD$-set for some integer $k\geq 1$. 
 In this paper, after presenting preliminaries, we  count the number of fair dominating sets of some specific graphs.  
\end{abstract}


\medskip
\noindent{\bf AMS Subj.\ Class.}: 05C25 

\section{Introduction and definitions}
Let $G = (V,E)$ be a simple graph with $n$ vertices. The distance between two vertices $u$ and $v$ denoted by $d(u,v)$ is the number of edges in a shortest path (also called a graph geodesic) connecting them.  
Let $S\subseteq  V$ be any subset of vertices of $G$. The induced subgraph $G[S]$ is the graph whose vertex set is $S$ and whose edge set consists of all of the edges in $E$ that have both endpoints in $S$. 
 A set $D\subseteq V(G)$ is a  dominating set, if every vertex in $V(G)\backslash D$ is adjacent to at least one vertex in $D$.
The  domination number $\gamma(G)$ is the minimum cardinality of a dominating set in $G$. 

For $k \geq 1$, a $k$-fair dominating set ($kFD$-set) in $G$, is a dominating set $D$ such that $|N(v) \cap D|=k$ for every vertex $ v \in V\setminus D$.
 The $k$-fair domination number of $G$,  denoted by $fd_k(G)$, is the minimum cardinality of a $kFD$-set. A $kFD$-set of $G$
 of cardinality $fd_k(G)$ is called a $fd_k(G)$-set. A fair dominating set, abbreviated
 FD-set, in $G$ is a $kFD$-set for some integer $k\geq 1$. The fair domination number,
 denoted by $fd(G)$, of a graph $G$ that is not the empty graph is the minimum
 cardinality of an FD-set in $G$. An FD-set of $G$ of cardinality $fd(G)$ is called a
 $fd(G)$-set. 
By convention, if $G=\overline{K_n}$, we define $fd(G)=n$. By the definition it is easy to see that for any graph $G$ of order $n$, $\gamma(G)\leq fd(G)\leq n$ and 
$fd(G)=n$ if and only if $G=\overline{K_n}$. 
Caro, Hansberg and Henning in \cite{Henning} showed that for a disconnected graph $G$ (without isolated vertices) of order $n\geq 3$, $fd(G)\leq n-2$, and they constructed 
an infinite family of graphs achieving equality in this bound.

The  corona of two graphs $G_1$ and $G_2$, is the graph
$G=G_1 \circ G_2$ formed from one copy of $G_1$ and $|V(G_1)|$ copies of $G_2$,
where the ith vertex of $G_1$ is adjacent to every vertex in the ith copy of $G_2$.
The corona $G\circ K_1$, in particular, is the graph constructed from a copy of $G$,
where for each vertex $v\in V(G)$, a new vertex $v'$ and a pendant edge $vv'$ are added.
The  join of two graphs $G_1$ and $G_2$, denoted by $G_1\vee G_2$,
is a graph with vertex set  $V(G_1)\cup V(G_2)$
and edge set $E(G_1)\cup E(G_2)\cup \{uv| u\in V(G_1)$ and $v\in V(G_2)\}$. Caro, Hansberg, and  Henning in \cite{Henning} proved that if $T$ is a tree of order $n\geq 2$, then 
$fd(T)\leq \frac{n}{2}$ with equality if and only if $T=T'\circ K_1$ for some tree $T'$. 

We know that if $S$ is a dominating set of $G$ and $S\subseteq S'$, then $S'$ is a dominating set, too. But this is not true for the fair dominating sets. As an example consider 
the cycle $C_9$ with $V(C_9)=\{1,2,...,9\}$. Observe that there are three fair dominating sets with cardinality three for $C_9$, but 
there is no dominating set of $C_9$ with cardinality four. This notation shows that study the fair dominating sets and finding the number of fair dominating sets of a graph with arbitrary cardinality is not easy problem. 

Regarding to enumerative side of dominating sets, Alikhani and Peng in  \cite{saeid1}, have introduced the domination polynomial of a graph. The domination polynomial of graph $G$ is the  generating function for the number of dominating sets of  $G$, i.e., $D(G,x)=\sum_{ i=1}^{|V(G)|} d(G,i) x^{i}$ (see \cite{euro,saeid1}).   This  polynomial and its roots has been actively studied in recent
years (see for example \cite{Brown,Kot}). 
It is natural to count the number of another kind of dominating sets (\cite{utilitas}).     
Let ${\cal D}_f(G,i)$ be the family of the 
fair   dominating sets of a graph $G$ with cardinality $i$ and let
$d_f(G,i)=|{\cal D}_f(G,i)|$. In this paper we count the number the fair dominating sets of certain graphs.  We denote the set $\{1,2,...,n\}$ simply by $[n]$.

\section{Main results}

In this section, similar to the domination polynomial, we state the definition of the fair domination polynomial and some
of its properties. Then, we count the number of fair dominating sets of specific graphs such as complete bipartite graph $K_{n,n}$, cycles, paths, friendship graph and triangular cactus chain. 

\subsection{Fair domination polynomial} 

In this subsection, we state the definition of the fair domination polynomial and some
of its properties.
\begin{definition}\label{Definition2}
	Let ${\cal D}_f(G,i)$ be the family of the fair dominating sets of a graph $G$ with cardinality $i$ and let
	$d_f(G,i)=|{\cal D}_f(G,i)|$. Then the fair domination polynomial $D_f(G,x)$ of $G$ is defined as
	\begin{center}
		$D_f(G,x)=\displaystyle\sum_{i=fd(G)}^{|V(G)|} d_f(G,i) x^{i}$,
	\end{center}
	where $fd(G)$ is the fair domination number of $G$.
\end{definition}

As an example, suppose that $V(K_n)=[n]$. Since every subsets of size $i\geq1$ is a fair dominating set of $K_n$, so we have  $D_f(K_n,x)=(1+x)^n-1$. Note that in the domination polynomial of a graph $G$,  every numbers in  $\{\gamma, \gamma+1,...,n\}$ appears as a power of $x$, but this is not true for the fair domination polynomial of any graph. For example in the fair domination polynomial of $C_9$, the term $x^6$ does not appear. As another example, in the $D_f(K_{4,4},x)$, the term $x^3$ does not appear, but for instance the term $x^2$ exists with coefficient  16.    The following theorem is easy and follows from the definition of the fair domination polynomial.    
\begin{theorem}\label{theorem42}
	Let $G$ be a graph with $|V(G)|=n$. Then
	\begin{enumerate}
		\item[(i)] If $G$ is connected, then $d_f(G,n)=1$ and $d_f(G,n-1)=n$,
		\item[(ii)] If $i<fd(G)$ or $i>n$, then $d_f(G,i)=0$.  
		\item[(iii)] $D_f(G,x)$ has no constant term.
		\item[(iv)] Zero is a root of $D_f(G,x)$, with multiplicity $fd(G)$.
	\end{enumerate}
\end{theorem}

\subsection{Results for $K_{n,n}$ and  $C_n$} 
In this subsection, we study the number of fair dominating sets of complete bipartite graph $K_{n,n}$ and the cycle graph $C_n$. We start with $K_{n,n}$. 

\begin{theorem}
	\begin{enumerate}  
		\item[(i)]  If $r>2$ is odd, then 
		\[
	d_f(K_{n,n},r)=\left\{
		\begin{array}{lr}
		{\displaystyle 2};&
		\quad\mbox{if $r=n$,}\\[15pt]
		{\displaystyle 0};&
		\quad\mbox{if $r<n$,}\\[15pt]
		{\displaystyle 2{n\choose r-n}};&
		\quad\mbox{if $r>n$.}
		\end{array}
		\right.
		\] 
		
		\item[(ii)]  If $r\geq 2$ is even, then 
		\[
		d_f(K_{n,n},r)=\left\{
		\begin{array}{lr}
		{\displaystyle {n\choose r/2}^2};&
			\quad\mbox{if $r<n$,}\\[15pt]
			{\displaystyle {n\choose r/2}^2+2};&
			\quad\mbox{if $r=n$,}\\[15pt]
			{\displaystyle {n\choose r/2}^2+{n\choose r-n}^2};&
				\quad\mbox{if $r>n$.}
				\end{array}
				\right.
				\] 
		\end{enumerate} 
		\end{theorem} 
\proof
\begin{enumerate}  
	\item[(i)]
	Let $X$ and $Y$ be two parts of $K_{n,n}$ and $D$ be a fair dominating set of $K_{n,n}$ with cardinality odd $r$. Let $D=D_1\cup D_2$, where $D_1\subseteq X$ and $D_2\subseteq Y$. Since $r$ is odd and $r=n$, we can suppose that $|D_1|$ is odd and $|D_2|$ is even and so the  vertices of $Y\setminus D_2$ are dominated by the odd number of vertices and vertices of $X\setminus D_1$  are dominated by the even number of vertices. Therefore, the  fair
	dominating set of cardinality $r$  is either entire $X$ or $Y$ and so we have the result for this case. 
	If $r<n$, with the same explanation, there is no fair dominating set of cardinality $r$ for $K_{n,n}$. Now, suppose that $r>n$. Choosing all vertices of part $X$ (or part $Y$) and $r-n$ vertices from $Y$ (or from $X$) gives a fair dominating set of size $r$. So if $r>n$, then $d_f(K_{n,n},r)= 2{n\choose r-n}$.  
	
	\item[(ii)]      
	First suppose that $r<n$. Choosing all $r$ vertices from $X$ or $Y$ is not a  fair dominating set of $K_{n,n}$. If $D=D_1\cup D_2$, where $D_1\subseteq X$ and $D_2\subseteq Y$, then both $|D_1|$ and $|D_2|$ are either odd or even, because $|D|=r$ is even. Obviously, if $|D_1|\neq |D_2|$, then $D$ cannot be a fair dominating set. So by choosing $\frac{r}{2}$ vertices from each part we have a fair dominating set. So the number of fair dominating sets of $K_{n,n}$ with size $r<n$ is ${n\choose \frac{r}{2}}$. 
	 For the case $r=n$, the parts $X$ and $Y$ can be fair dominating sets, too. So we have the result for case $r=n$. 
	 For the case $r>n$, every fair dominating sets $D$ of $K_{n,n}$ can be constructed  by choosing $\frac{r}{2}$ vertices from each part and also by choosing $n$ vertices from one part and $r-n$ vertices from another part. Therefore the results follows. \qed  
	 
	\end{enumerate}

Here we consider the number of fair dominating sets of cycles.  
 Let $C_n, n\geq 3$, be the cycle with $n$ vertices $V(C_n)=[n]$ and $E(C_n)=\{\{1,2\},\{2,3\},...,\{n-1,n\},\{n,1\}\}$.  Let $f\mathcal{C}_n^i$ be the family of fair dominating sets of $C_n$ with cardinality $i$. We shall investigate the fair dominating sets of cycles.
 A {\it simple path} is a path where all its internal vertices have degree two. We need the following easy lemma to prove our main results in this section:

\begin{lemma}\label{lemma1}
 The following properties hold for cycles,
 	
 	\begin{enumerate}
 	
 		\item[(i)]
 		{\rm(\cite{Henning})} 		
 		$fd(C_{n})=\gamma(C_n)=\lceil\frac{n}{3}\rceil$, unless $n\equiv 2$ $(mod\, 3)$ and 
 		$n\geq 5$ in which case $fd(C_n)=\gamma(C_n)+1=\lceil\frac{n}{3}\rceil+1$.
 		
 		 		\item[(ii)]
 		$f\mathcal{C}_{j}^{i}=\emptyset$, if and only if $i>j$ or $i<\lceil\frac{j}{3}\rceil$.(by $(i)$ above).
 		
 		\item[(iii)]
 	If a graph $G$ contains a simple path of length $3k-1$, then every fair dominating set of $G$ must contain at least $k$ vertices of the path.
 		
 		\end{enumerate}
 	\end{lemma}

To construct the fair dominating sets of $C_n$ with $V(C_n)=[n]$ of size $k$, we partition the number $k$ to $n-k$ natural numbers, when $n-k$ is odd or $n-k=k$ and  partition the number $k$ to $\frac{n-k}{2}$ natural numbers, when $n-k$ is even. 
Suppose that $k=t_1+t_2+...+t_{n-k}$, where $t_i$ ($1\leq i\leq n-k$) is a positive   integer number.
We define the family $\mathcal{A}\subseteq [n]$, based on partition  $k=t_1+t_2+...+t_{n-k}$ as follows: 

	$$\mathcal{A}=A_1\cup A_2\cup...\cup A_{n-k},$$
		where the set $A_1$ contains $t_1$ consecutive numbers from $[n]$, the set $A_2$ contains  $t_2$ consecutive numbers from $[n]\setminus A_1$ and finally the set $A_{n-k}$ contains  $t_{n-k}$ consecutive numbers from $[n]\setminus \big(A_1\cup A_2\cup...\cup A_{n-k-1}\big)$ and also $d(G[A_i],G[A_{i+1}])=2$.  
	 
\medskip	
	Now suppose that $k=t_1+t_2+...+t_{\frac{n-k}{2}}$, where $t_i$ ($1\leq i\leq \frac{n-k}{2}$) is a non-negative integer number. We define the family $\mathcal{B}\subseteq [n]$, based on this 
	 partition, as follows: 
	
	$$\mathcal{B}=B_1\cup B_2\cup...\cup B_{\frac{n-k}{2}},$$
	where the set $B_1$ contains $t_1$ consecutive numbers from $[n]$, the set $B_2$ contains  $t_2$ consecutive numbers from $[n]\setminus B_1$ and finally the set $B_{\frac{n-k}{2}}$ contains  $t_{\frac{n-k}{2}}$ consecutive numbers from $[n]\setminus \big(B_1\cup B_2\cup...\cup B_{\frac{n-k}{2}-1}\big)$ and also $d(G[B_i],G[B_{i+1}])=3$.

	\medskip 
 
 Now with these constructions and explanations, the family of fair dominating sets of $C_n$ follows from   the following theorem. 

\begin{theorem}  \label{1}
	Let $C_n, n\geq 3$, be the cycle with $n$ vertices. The following parts give the fair dominating sets of $C_n$ with cardinality $k$.  
	\begin{enumerate}
		\item [(i)] If $n-k$ is even and $n\leq 2k$, then  $f\mathcal{C}_n^k=\mathcal{A}\cup \mathcal{B}$.

		\item[(ii)] 
		If $n-k$ is even and $n>2k$, then  $f\mathcal{C}_n^k=\mathcal{B}$.
		
		\item[(iii)]  
		If $n-k$ is odd and $n\leq 2k$, then  $f\mathcal{C}_n^k=\mathcal{A}$. 	
		
		\item[(iv)]  
		If $n-k$ is odd and $n>2k$, then  $f\mathcal{C}_n^k=\emptyset$.

		\end{enumerate} 
\end{theorem}

Let explain Theorem \ref{1} with the following example:

\begin{example} 
Consider the cycle $C_8$ with $V(C_8)=[8]$. We want to obtain $d_f(C_8,4)$. Here $n=8$ and $k=4$ and so $n-k=4$. Since $n-k$ is even and $n-k=k$, so by Theorem \ref{1}, we have 
$f\mathcal{C}_8^4=\mathcal{A}\cup \mathcal{B}$, where $\mathcal{A}$ and $\mathcal{B}$ constructed as follows. 

\noindent First partition $n-k=4$ to $k=4$ positive integer umbers, i.e., $4=1+1+1+1$. So by explanation before Theorem \ref{1}, $\mathcal{A}=\{1,3,5,7\}\cup\{2,4,6,8\}$. 
Now we partition $n-k=4$ to $\frac{n-k}{2}=2$ positive integers, i.e., $4=1+3$ and $4=2+2$. 
For the partition $4=1+3$, we have the following fair dominating sets 
$$\{i\}\cup\{i+3,i+4,i+5\},$$
  where $i\in [8]$, and so we have eight fair dominating sets in this case. 
  
 \noindent Now we consider the partition $4=2+2$. From this partition we have the following fair dominating sets: 
   $$\{i,i+1\}\cup\{i+4,i+5\},$$
   where $i\in [8]$, and so we have four  fair dominating sets in this case. Therefore  $d_f(C_8,4)=2+8+4=14$. 
	\end{example}

Here, by Theorem \ref{1}, we count the number of the fair dominating sets of cycles. 
Constructing the fair dominating sets of cycle $C_n$ with cardinality $k$ is similar to the following known problem. Assume that we have a set $X$ with $n-k$ elements and we want to choose a subset $A_1$ of $X$ containing $t_1$ elements, then choosing a subset $A_2$ from $X\setminus A_1$ containing  $t_2$ elements and finally  choosing a subset $A_{n-k}$ from $X\setminus (A_1\cup...\cup A_{n-k-1})$ containing  $t_{n-k}$ elements. How many cases we have? As we know this is the generalized combination and this number is denoted by 
${n-k\choose t_1,t_2,...,t_{n-k}}$ which is equal to $\frac{(n-k)!}{t_1!t_2!...t_{n-k}!}$, where $n-k=t_1+t_2+...+t_{n-k}$. Note that since the elements in our study are in the around of  the roundtable, the number is equal to $\frac{1}{n-k}{n-k\choose t_1,t_2,...,t_{n-k}}$. 
In summary, if $P_j$  is a partition of the number $k$ to $n-k$ numbers, then 
$$|\mathcal{A}|=\sum_{P_j} \frac{1}{n-k}{n-k\choose t_1,t_2,...,t_{n-k}}.$$

With similar statements, 
$$|\mathcal{B}|=\sum_{P_j} \frac{2}{n-k}{\frac{n-k}{2}\choose t_1,t_2,...,t_{\frac{n-k}{2}}}.$$
Now, by Theorem \ref{1} we state the following results:

\begin{theorem}  \label{2}
	Let $C_n, n\geq 3$, be the cycle of order  $n$.   
	\begin{enumerate}
		\item [(i)] If $n-k$ is even and $n\leq 2k$, then  $d_f(C_n,k)=n(|\mathcal{A}|+|\mathcal{B}|)$.

		\item[(ii)] 
		If $n-k$ is even and $n>2k$, then  $d_f(C_n,k)=n|\mathcal{B}|$.
		
		\item[(iii)]  
		If $n-k$ is odd and $n\leq 2k$, then  $d_f(C_n,k)=n|\mathcal{A}|$. 	
		
		\item[(iv)]  
		If $n-k$ is odd and $n>2k$, then  $d_f(C_n,k)=0$. 
			
		\end{enumerate} 
	\end{theorem} 
	
	Using Theorem \ref{2}, we obtain $d_f(C_n,j)$ for $1\leq n\leq 12$ as shown in Table 1. Using this table we have the following corollary. 
	
	\begin{corollary} 
		\begin{enumerate} 
			\item[(i)] For every $n\geq 3$, $d_f(C_n,n-2)=\frac{(n-1)n}{2}$.
			\item [(ii)]For every $n\geq 6$, $d_f(C_n,n-3)=\frac{(n-5)(n-4)n}{6}$.
			\item[(iii)] 
			For every $n\geq 7$, $d_f(C_n,fd(C_n))=\frac{(n-8)(n-7)n}{6}$.   
			
			 \end{enumerate}
		\end{corollary} 
			
\[
		\begin{footnotesize}
		\small{
			\begin{tabular}{r|lcrrrcccccccccccc}
			$j$&$1$&$2$&$3$&$4$&$5$&$6$&$7$&$8$&$9$&$10$&$11$&$12$\\[0.3ex]
			\hline
			$n$&&&&&&&&&&&&&&&&\\
			$3$&3&3&1&&&&&&&&&&&&&\\
			$4$&0&6&4&1&&&&&&&&&&&&\\
			$5$&0&0&10&5&1&&&&&&&&&&&\\
			$6$&0&3&2&9&6&1&&&&&&&&&&\\
			$7$&0&0&7&7&21&7&1&&&&&&&&&\\
			$8$&0&0&0&14&16&28&8&1&&&&&&&&\\
			$9$&0&0&3&0&27&30&36&9&1&&&&&&&\\
			$10$&0&0&0&10&2&50&50&45&10&1&&&&&&\\
			$11$&0&0&0&0&22&11&88&77&55&11&1&&&&&\\
			$12$&0&0&0&3&0&42&36&147&112&66&12&1&&&&\\
			\end{tabular}}
		\end{footnotesize}
		\]
		\begin{center}
			\noindent{Table 1.} $d_f(C_{n},j)$, the number of fair dominating sets of $C_n$ with cardinality $j$.
		\end{center}
		
\[
		\begin{footnotesize}
		\small{
			\begin{tabular}{r|lcrrrcccccccccccc}
			$j$&$1$&$2$&$3$&$4$&$5$&$6$&$7$&$8$&$9$&$10$&$11$&$12$\\[0.3ex]
			\hline
			$n$&&&&&&&&&&&&&&&&\\
			$1$&1&&&&&&&&&&&&&&&\\
			$2$&2&1&&&&&&&&&&&&&&\\
			$3$&1&3&1&&&&&&&&&&&&&\\
			$4$&0&2&4&1&&&&&&&&&&&&\\
			$5$&0&2&4&5&1&&&&&&&&&&&\\
			$6$&0&1&4&7&6&1&&&&&&&&&&\\
			$7$&0&0&3&7&11&7&1&&&&&&&&&\\
			$8$&0&0&2&6&12&16&8&1&&&&&&&&\\
			$9$&0&0&1&3&11&20&22&9&1&&&&&&&\\
			$10$&0&0&0&3&12&20&32&29&10&1&&&&&&\\
			$11$&0&0&0&2&6&21&36&49&37&11&1&&&&&\\
			$12$&0&0&0&1&8&20&36&64&60&46&12&1&&&&\\
			\end{tabular}}
		\end{footnotesize}
		\]
		\begin{center}
		\label{2}	\noindent{Table 2.} $d_f(P_{n},j)$, the number of fair dominating sets of $P_n$ with cardinality $j$.
		\end{center}		

\subsection{Results for $P_n$, friendship graph $F_n$ and triangular cactus $T_n$} 
  
Using Maple, we computed the number of fair dominating sets of path graph $P_n$ for $n\leq 12$ and have shown them in Table 2. From this table we have the following 
easy results. Unfortunately, our attempt for finding a closed formula for the number of fair dominating sets of $P_n$ failed by now.  

\begin{theorem} For every $n\geq 2$, 
	\begin{enumerate} 
\item[(i)] $d_f(P_n,n)=1,$~$d_f(P_{n},n-1)=n$. 
\item[(ii)] $d_f(P_{n},n-2)=\frac{(n-3)(n-2)}{2}+1$.
\item[(iii)] $d_f(P_{n},n-4)=\displaystyle\sum_{i=1}^{n-3} {n-3-i \choose n-3-3i}.$
\end{enumerate} 
\end{theorem}	

\begin{theorem} For every $k\geq 2$,
	\begin{enumerate} 
		\item[(i)] $d_f(P_{3k},k)=1$. 
		\item[(ii)] $d_f(P_{3k+1},k+1)=3$.
		\item[(iii)] $d_f(P_{3k+2},k+1)=2.$
	\end{enumerate} 
\end{theorem}


Here, we consider the number of fair dominating sets of friendship graphs. The friendship graph $F_n$,  is the graph obtained by taking $n$ copies of the cycle graph $C_3$ with a vertex in common. 
It is easy to see that we have no fair dominating set of even size for $F_n$. So we have the following result: 

\begin{theorem}
For every $1\leq i\leq \lfloor\frac{n}{2}\rfloor$, $d_{f}(F_n,2i)=0.$
\end{theorem}

The following theorem gives the number of fair dominating sets of $F_n$ with odd cardinality. 

\begin{theorem}
	\begin{enumerate} 
\item[(i)]
For all $n\geq 3$, $d_{f}(F_n,3)=n.$
\item[(ii)] For every $n\geq 5$, 
$d_{f}(F_n,5)=\frac{n(n-1)}{2}.$
\end{enumerate} 
\end{theorem} 
\proof
\begin{enumerate}
\item[(i)]
The only way to make $FD$-sets of $F_n$ of size $3$ is to choose all vertices of one triangle, whereas there are $n$ triangles, the number of $FD$-sets of size three is  $n$.
\item[(ii)]
The ways to construct $FD$-sets of graph $F_n$ of size $5$ are as follows;\\
 We choose the first $FD$-set in the way that includes all vertices of triangles $1$ and $2$,
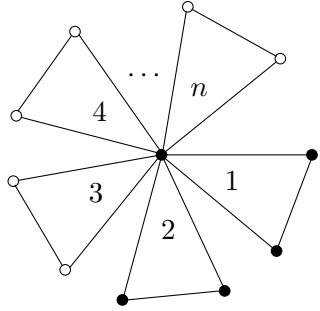
\begin{figure} 
\begin{center}
\begin{tikzpicture}
\filldraw (0,0) circle (2pt);
\filldraw(2,0) circle(2pt);
\draw(1.58,1.28) circle(2pt);
\draw(0.35,1.97) circle(2pt);
\draw(-1.15,1.64) circle(2pt);
\draw(-1.93,0.52) circle(2pt);
\draw(-1.97,-0.35) circle(2pt);
\draw(-1.28,-1.53) circle(2pt);
\filldraw(-0.52,-1.93) circle(2pt);
\filldraw(0.84,-1.81) circle(2pt);
\filldraw(1.53,-1.28) circle(2pt);

\draw(1.94,0)--(0.06,0);
\draw(1.53,1.24)--(0.05,0.04);
\draw(0.33,1.91)--(0.02,0.06);
\draw(-1.1,1.6)--(-0.02,0.06);
\draw(-1.87,0.5)--(-0.06,0.03);
\draw(-1.91,-0.34)--(-0.06,-0.01);
\draw(-1.22,-1.5)--(-0.05,-0.05);
\draw(-0.5,-1.87)--(-0.01,-0.06);
\draw(0.82,-1.75)--(0.03,-0.06);
\draw(1.48,-1.24)--(0.05,-0.05);

\draw(1.99,-0.06)--(1.55,-1.22);
\draw(1.53,1.31)--(0.4,1.94);

\draw(-1.2,1.58)--(-1.91,0.58);

\draw(-1.95,-0.4)--(-1.32,-1.47);

\draw(-0.46,-1.93)--(0.78,-1.81);

\draw(0.5, 0.87)circle (0pt) node{$n$};

\draw(0.94, -0.34)circle (0pt) node{$1$};

\draw(0.09,-1)circle (0pt) node{$2$};
\draw(-0.87,-0.5)circle (0pt) node{$3$};
\draw(-0.82,0.57)circle (0pt) node{$4$};

\draw(-0.2,1.05)circle (0pt) node{$\cdots$};
\end{tikzpicture}
\end{center}
\caption{\label{size5} A fair dominating set of $F_n$ of cardinality $5$.}
\end{figure} 
 then by fixing triangle $1$ and shifting forward triangle $2$, we have another $FD$-set (see Figure \ref{size5}). By continue this process there will be $n-1$ cases (for example triangles $1$ and $3$, triangles $1$ and $4$,$\cdots$). In the next step by one shifting forward triangles $1$ and $2$ we have a new $FD$-set. Also by fixing triangle $2$ and shifting forward triangle $3$, a new $FD$-set is made. By continuing  the process there are $n-2$ cases. We repeat this action until we construct the last $FD$-set. So the number of fair dominating sets of $F_n$ of size $5$ is
$$(n-1)+(n-2)+(n-3)+\cdots+3+2+1,$$
which equals ${n(n-1) \over 2}.$
\end{enumerate}\qed

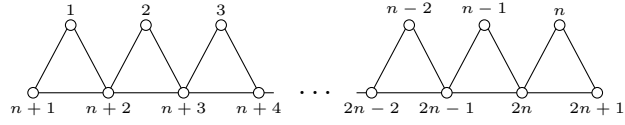
\begin{figure} 
 \begin{center}
\begin{tikzpicture}
\draw (0,0) circle (2pt) node[below]{\tiny{${n+1}$}};
\draw (1,0) circle (2pt) node[below]{\tiny{${n+2}$}};
\draw (2,0) circle (2pt) node[below]{\tiny{${n+3}$}};
\draw (3,0) circle (2pt) node[below]{\tiny{${n+4}$}};
\draw (4.5,0) circle (2pt) node[below]{\tiny{${2n-2}$}};
\draw (5.5,0) circle (2pt) node[below]{\tiny{${2n-1}$}};
\draw (6.5,0) circle (2pt)node[below]{\tiny{${2n}$}};
\draw (7.5,0) circle (2pt) node[below]{\tiny{${2n+1}$}};
\draw (0.5,0.9) circle (2pt) node[above]{\tiny{${1}$}};
\draw (1.5,0.9) circle (2pt) node[above]{\tiny{${2}$}};
\draw (2.5,0.9) circle (2pt) node[above]{\tiny{${3}$}};
\draw (7,0.9) circle (2pt) node[above]{\tiny{${n}$}};
\draw (5,0.9) circle (2pt) node[above]{\tiny{${n-2}$}};
\draw (6,0.9) circle (2pt) node[above]{\tiny{${n-1}$}};
\draw (0.07,0)--(0.93,0);
\draw (1.07,0)--(1.93,0);
\draw (2.07,0)--(2.93,0);
\draw (3.07,0)--(3.2,0);

\draw (3.75,0) circle (0pt) node{$\ldots$};
\draw (0.03,0.06)--(0.45,0.85);
\draw (0.55,0.85)--(0.97,0.06);
\draw(1.03,0.06)--(1.45,0.85);
\draw (1.55,0.85)--(1.97,0.06);
\draw(2.03,0.06)--(2.45,0.85);
\draw (2.55,0.85)--(2.97,0.06);
\draw(4.53,0.06)--(4.95,0.85);
\draw (5.06,0.85)--(5.47,0.06);
\draw(5.53,0.06)--(5.95,0.85);
\draw (6.05,0.85)--(6.47,0.06);
\draw(6.53,0.06)--(6.95,0.85);
\draw (7.05,0.85)--(7.47,0.06);
\draw (4.3,0)--(4.43,0);
\draw(4.57,0)--(5.43,0);
\draw(5.57,0)--(6.43,0);
\draw(6.57,0)--(7.43,0);
\end{tikzpicture}
\end{center}
\caption{\label{T_n} The graph $T_n$.  }
\end{figure}

Cactus graphs, were first known as Husimi tree, they appeared in the scientific literature some sixty years ago in papers by Husimi and
Riddell concerned with cluster integrals in the theory of condensation in statistical mechanics \cite{9,12,14}.   Cactus graph is a connected graph in which any two graph cycles have no edge in common. Equivalently, it is a connected graph in which any two (simple) cycles have at most one vertex in common.

The following theorem gives the number of fair dominating sets of triangular cactus graph $T_n$ of some specific sizes (see Figure \ref{T_n}).

\begin{theorem}
\begin{enumerate}
\item[(i)]
 $d_f(T_n,2n+1)=1$.
\item[(ii)]
$d_f(T_n,2n)= {2n+1 \choose 1}$.

\item[(iii)]
$d_f(T_n,2n-1)={n \choose {n-2}}+{{n-1} \choose {n-3}}+2n.$

\item[(iv)]
$d_f(T_n,2n-2)={ n \choose {n-3}}+2 {{n-1} \choose {n-3}}.$
\item[(v)]
$d_f(T_n, 2n-3)={n \choose n-4}+1+d_f(P_{n+1},n-3)~~~~~(n \neq 4).$
 
  \item[(vi)]
  	$
  	d_f(T_n,2n-4)=\left\{
  	\begin{array}{lr}
  	{\displaystyle {n \choose n-5}};&
  	\quad\mbox{if $6 \leq n \leq 9$,}\\[15pt]
    	{\displaystyle {n \choose n-5}+d_f(P_{n+1},n-4)};&
  	\quad\mbox{if $n > 9$.}
  	\end{array}
  	\right.
  	$ 
   \item[(vii)]
$d_f(T_n, n+1)=4$, ~~$(n \neq 1, 2, 4)$.
\item[(viii)]
$d_f(T_n,n)=0$,~~ $(n \neq 1,3)$.
\end{enumerate}
\end{theorem}
\proof
\begin{enumerate}
	\item[(i)] 
	The proof is  straightforward. 
	\item[(ii)] 
			The proof is  straightforward.
	\item[(iii)]
We can construct $FD$-sets of $T_n$ of size $2n-1$ in three ways. In the first method, we choose $n+1$ vertices of down vertices  of $T_n$ (i.e., vertices of induced path $P_{n+1}$) and then to reach size $2n-1$, we have ${ n \choose {n-2}}$ choices from the  top  vertices of $T_n$.

 In the second method for constructing $FD$-sets of $T_n$,  we choose vertices as follows: We choose all $n$ vertices in the top of cactus $T_n$ (i.e., all vertices in $\{1,2,...,n\}$), and  the first and the last vertices of induced path $P_{n+1}$ in $T_n$.  
 Then to reach the size of $2n-1$, we have $ {{n-1} \choose {n-3}}$ choices of down  vertices of $T_n$.  
 In the third method,  choose vertices as follows. To construct $\mathcal{D}_f(T_n,2n-1)$, it is suffices to remove two vertices of $V(T_n)$. These two vertices should choose from the following set: 
   $$\lbrace (n+1,n),(n+1, n-1), (n+1, n-2), (n+1, n-3), \cdots, (n+1,2), (n+1, 1)\rbrace,$$
  and so we have $n$, $FD$-sets by these unchosen vertices. 
 It is easy to see that we can choose two vertices from the following set, too. 
   $$\lbrace (1, 2n+1), (2, 2n+1), (3, 2n+1), \cdots, (n-2, 2n+1), (n-1, 2n+1), (n, 2n+1) \rbrace.$$
Therefore, by additive principle we have $d_f(T_n,2n-1)={n \choose {n-2}}+{{n-1} \choose {n-3}}+2n.$

\item[(iv)]
To make $FD$-sets of size $2n-2$, first we choose $n+1$ vertices of down vertices  of $T_n$ (i.e., all vertices of $\{n+1,n+2,....,2n,2n+1\}$), then to reach size $2n-2$ we need $n-3$ vertices, that we can choose them by $ {n \choose {n-3}}$ choices.

In the second  way, first we choose $n+1$ vertices as follows. Choose the first vertex (or the last vertex) from the top of $T_n$ and $n$ vertices of the down vertices of $T_n$ except the first vertex (last vertex) of down. In other words we choose vertices $\{1,n+2,n+3,...,2n, 2n+1\}$ or $\{n,n+1,n+2,...,2n-1,2n\}$. 
 So we have the result.

\item[(v)]
To construct fair dominating sets of $T_n$ of size $2n-3$, we have three  cases:\\
\textbf{Case 1:} We choose all vertices of down vertices of $T_n$, i.e.,  $\lbrace {n+1},{n+2}, \cdots,2n+1 \rbrace$, then we have ${ n \choose n-4}$ choices. 
 
 \textbf{Case 2:} The second  way is choosing all vertices from the top of $T_n$, i.e, $\{1, 2,..., n-1, n\}$. 
 Now we shall choose  $n-3$ vertices of induced path $P_{n+1}$ that we have seen that the number of choices is  $\sum _{P_j}{n-3\choose 3,5}$.
 
  \textbf{Case 3:} It is suffices to do not choose the four vertice $\{1,n+1,n,2n+1\}$. Therefore, we have the result.

\item[(vi)]
To construct FD-sets of size $2n-4$, one way is to choose $n+1$ down vertices of $T_n$, then we have ${n \choose {n-5}}$ choices to make FD-sets.

 Another way is choosing $n$ vertices from the top vertices of $T_n$ and choosing $n-4$ vertices of induced path $P_{n+1}$ which we have seen that its number equal  $\sum _{P_j}{n-4\choose 6}$ for $n>9$.

\item[(vii)]
To construct fair dominating sets in $T_n$ of size $n+1$ we have four ways. The first, is choosing all down vertices in $T_n$.   Another methods are the following sets (see Figure \ref{n+1}).   

$\{1,n,n+2,n+3,...,2n\}$, $\{1,n+2,n+3,...,2n,2n+1\}$, $\{n,n+1,...,2n\}$. 

So, we have the result.\qed

\begin{figure} 
\begin{center}
\begin{tikzpicture}
\draw (0,0) circle (2pt) node[below]{\tiny{${n+1}$}};
\filldraw (1,0) circle (2pt) node[below]{\tiny{${n+2}$}};
\filldraw (2,0) circle (2pt) node[below]{\tiny{${n+3}$}};
\filldraw (3,0) circle (2pt) node[below]{\tiny{${n+4}$}};
\filldraw (4.5,0) circle (2pt) node[below]{\tiny{${2n-2}$}};
\filldraw (5.5,0) circle (2pt) node[below]{\tiny{${2n-1}$}};
\filldraw (6.5,0) circle (2pt)node[below]{\tiny{${2n}$}};
\draw (7.5,0) circle (2pt) node[below]{\tiny{${2n+1}$}};
\filldraw (0.5,0.9) circle (2pt) node[above]{\tiny{${1}$}};
\draw (1.5,0.9) circle (2pt) node[above]{\tiny{${2}$}};
\draw (2.5,0.9) circle (2pt) node[above]{\tiny{${3}$}};
\filldraw (7,0.9) circle (2pt) node[above]{\tiny{${n}$}};
\draw (5,0.9) circle (2pt) node[above]{\tiny{${n-2}$}};
\draw (6,0.9) circle (2pt) node[above]{\tiny{${n-1}$}};
\draw (0.07,0)--(0.93,0);
\draw (1.07,0)--(1.93,0);
\draw (2.07,0)--(2.93,0);
\draw (3.07,0)--(3.2,0);

\draw (3.75,0) circle (0pt) node{$\ldots$};
\draw (0.03,0.06)--(0.45,0.85);
\draw (0.55,0.85)--(0.97,0.06);
\draw(1.03,0.06)--(1.45,0.85);
\draw (1.55,0.85)--(1.97,0.06);
\draw(2.03,0.06)--(2.45,0.85);
\draw (2.55,0.85)--(2.97,0.06);
\draw(4.53,0.06)--(4.95,0.85);
\draw (5.06,0.85)--(5.47,0.06);
\draw(5.53,0.06)--(5.95,0.85);
\draw (6.05,0.85)--(6.47,0.06);
\draw(6.53,0.06)--(6.95,0.85);
\draw (7.05,0.85)--(7.47,0.06);
\draw (4.3,0)--(4.43,0);
\draw(4.57,0)--(5.43,0);
\draw(5.57,0)--(6.43,0);
\draw(6.57,0)--(7.43,0);
\end{tikzpicture}
\end{center} 
 
 \begin{center}
\begin{tikzpicture}
\draw (0,0) circle (2pt) node[below]{\tiny{${n+1}$}};
\filldraw (1,0) circle (2pt) node[below]{\tiny{${n+2}$}};
\filldraw (2,0) circle (2pt) node[below]{\tiny{${n+3}$}};
\filldraw (3,0) circle (2pt) node[below]{\tiny{${n+4}$}};
\filldraw (4.5,0) circle (2pt) node[below]{\tiny{${2n-2}$}};
\filldraw (5.5,0) circle (2pt) node[below]{\tiny{${2n-1}$}};
\filldraw (6.5,0) circle (2pt)node[below]{\tiny{${2n}$}};
\filldraw (7.5,0) circle (2pt) node[below]{\tiny{${2n+1}$}};
\filldraw (0.5,0.9) circle (2pt) node[above]{\tiny{${1}$}};
\draw (1.5,0.9) circle (2pt) node[above]{\tiny{${2}$}};
\draw (2.5,0.9) circle (2pt) node[above]{\tiny{${3}$}};
\draw (7,0.9) circle (2pt) node[above]{\tiny{${n}$}};
\draw (5,0.9) circle (2pt) node[above]{\tiny{${n-2}$}};
\draw (6,0.9) circle (2pt) node[above]{\tiny{${n-1}$}};
\draw (0.07,0)--(0.93,0);
\draw (1.07,0)--(1.93,0);
\draw (2.07,0)--(2.93,0);
\draw (3.07,0)--(3.2,0);

\draw (3.75,0) circle (0pt) node{$\ldots$};
\draw (0.03,0.06)--(0.45,0.85);
\draw (0.55,0.85)--(0.97,0.06);
\draw(1.03,0.06)--(1.45,0.85);
\draw (1.55,0.85)--(1.97,0.06);
\draw(2.03,0.06)--(2.45,0.85);
\draw (2.55,0.85)--(2.97,0.06);
\draw(4.53,0.06)--(4.95,0.85);
\draw (5.06,0.85)--(5.47,0.06);
\draw(5.53,0.06)--(5.95,0.85);
\draw (6.05,0.85)--(6.47,0.06);
\draw(6.53,0.06)--(6.95,0.85);
\draw (7.05,0.85)--(7.47,0.06);
\draw (4.3,0)--(4.43,0);
\draw(4.57,0)--(5.43,0);
\draw(5.57,0)--(6.43,0);
\draw(6.57,0)--(7.43,0);
\end{tikzpicture}
\end{center}

\begin{center}
\begin{tikzpicture}
\filldraw (0,0) circle (2pt) node[below]{\tiny{${n+1}$}};
\filldraw (1,0) circle (2pt) node[below]{\tiny{${n+2}$}};
\filldraw (2,0) circle (2pt) node[below]{\tiny{${n+3}$}};
\filldraw (3,0) circle (2pt) node[below]{\tiny{${n+4}$}};
\filldraw (4.5,0) circle (2pt) node[below]{\tiny{${2n-2}$}};
\filldraw (5.5,0) circle (2pt) node[below]{\tiny{${2n-1}$}};
\filldraw (6.5,0) circle (2pt)node[below]{\tiny{${2n}$}};
\draw (7.5,0) circle (2pt) node[below]{\tiny{${2n+1}$}};
\draw (0.5,0.9) circle (2pt) node[above]{\tiny{${1}$}};
\draw (1.5,0.9) circle (2pt) node[above]{\tiny{${2}$}};
\draw (2.5,0.9) circle (2pt) node[above]{\tiny{${3}$}};
\filldraw (7,0.9) circle (2pt) node[above]{\tiny{${n}$}};
\draw (5,0.9) circle (2pt) node[above]{\tiny{${n-2}$}};
\draw (6,0.9) circle (2pt) node[above]{\tiny{${n-1}$}};
\draw (0.07,0)--(0.93,0);
\draw (1.07,0)--(1.93,0);
\draw (2.07,0)--(2.93,0);
\draw (3.07,0)--(3.2,0);

\draw (3.75,0) circle (0pt) node{$\ldots$};
\draw (0.03,0.06)--(0.45,0.85);
\draw (0.55,0.85)--(0.97,0.06);
\draw(1.03,0.06)--(1.45,0.85);
\draw (1.55,0.85)--(1.97,0.06);
\draw(2.03,0.06)--(2.45,0.85);
\draw (2.55,0.85)--(2.97,0.06);
\draw(4.53,0.06)--(4.95,0.85);
\draw (5.06,0.85)--(5.47,0.06);
\draw(5.53,0.06)--(5.95,0.85);
\draw (6.05,0.85)--(6.47,0.06);
\draw(6.53,0.06)--(6.95,0.85);
\draw (7.05,0.85)--(7.47,0.06);
\draw (4.3,0)--(4.43,0);
\draw(4.57,0)--(5.43,0);
\draw(5.57,0)--(6.43,0);
\draw(6.57,0)--(7.43,0);
\end{tikzpicture}
\end{center}
\caption{\label{n+1} Fair dominating sets of $T_n$ with cardinality $n+1$.   }
\end{figure}
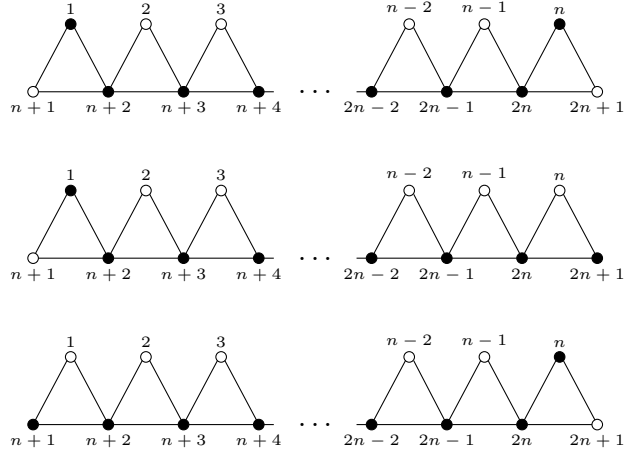

\end{enumerate}

\section{Conclusion} 

In this paper we studied the number of fair dominating sets of graphs. We obtained results for the complete bipartite graph $K_{n,n}$ and the cycle graph  $C_n$. But for path graph $P_n$, friendship graph $F_n$ and cactus triangle $T_n$ we have results for the number of fair dominating sets of specific size. There are many open problems in the study of the fair domination polynomial of a graph that we state and close the paper with some of them: 

\medskip
\noindent {\bf Problem 1}: Is there any recurrence relation for $d_f(G,k)$ based on the number of $FD$-sets of its subgraph, such as $G-v$, $G-e$, $G\circ e$ and ....?

\medskip
\noindent{\bf Problem 2}: Find a formula for the number of fair dominating sets of  path graph $P_n$, friendship graph $F_n$ and cactus triangle $T_n$.

\noindent{\bf Problem 3}: Characterize edges $e$ of graph $G$ such that $D_f(G,x)=D_f(G-e,x)$.


\end{document}